\documentclass[a4paper, authoryear]{article}
\usepackage{techRepsPDF}
\usepackage[parfill]{parskip}
\usepackage{natbib}
\usepackage[ruled,vlined,linesnumbered]{algorithm2e}
\SetCommentSty{textsf}
\usepackage{setspace} 
\usepackage{graphicx} 
\usepackage{placeins}
\usepackage{newproof}
\usepackage{setspace}
\usepackage{amsmath, amssymb}
\usepackage{xcolor, colortbl, color}
\usepackage{threeparttable, longtable}
\usepackage[margin=1cm]{caption}
\usepackage[pdftitle={Cutting Stock with Binary Patterns: Arc-flow Formulation}, pdfauthor={Filipe Brandao}]{hyperref}
\usepackage{tikz}
\usepackage{subfig}
\usetikzlibrary{shapes.geometric}
\usetikzlibrary{shapes.multipart}
\usetikzlibrary{decorations,decorations.markings}
\usetikzlibrary{shapes,arrows,fit,calc,positioning}
\usepackage{fullpage}
\usepackage{comment}

\newcolumntype{$}{>{\global\let\currentrowstyle\relax}} 
\newcolumntype{^}{>{\currentrowstyle}}

\def\HD#1#2{\vrule height #1pt depth #2pt width 0pt\relax}
\def\up{\HD{10}{0}}
\def\down{\HD{0}{5}}

\newcommand{\slp}{\mathop{\mathrm{lp}}}
\newcommand{\sip}{\mathop{\mathrm{ip}}}
\newcommand{\tip}{t^{\saf}_{\sip}}

\newcommand{\ncols}{\mathop{\#\mathrm{cols}}}

\newcommand{\saf}{\mathop{\mathrm{af}}}
\newcommand{\zlp}{z_{\slp}}
\newcommand{\zip}{z_{\sip}}
\newcommand{\tlp}{t^{\saf}_{\slp}}

\newcommand{\sgg}{\mathop{\mathrm{gg}}}
\newcommand{\tgg}{t^{\sgg}}

\newcommand{\bestlp}{\#{\saf}}
\newcommand{\bestgg}{\#{\sgg}}

\newcommand{\BigO}{\mathcal{O}}
\newcommand{\vS}{\textsc{s}}
\newcommand{\vT}{\textsc{t}}

\def\ptitleA{
Cutting Stock with Binary Patterns:\\
Arc-flow Formulation with\\ Graph Compression
}
\def\ptitleB{
Cutting Stock with Binary Patterns:\\
Arc-flow Formulation with Graph Compression
}
\def\pnumber{DCC-2013-09}
\def\pauthor{
Filipe Brand\~ao\\
{\small INESC TEC and Faculdade de Ci\^{e}ncias,
Universidade do Porto, Portugal}\\
{\small\texttt{fdabrandao@dcc.fc.up.pt}}
\ \\ 
\ \\
Jo\~ao Pedro Pedroso\\
{\small INESC TEC and Faculdade de Ci\^{e}ncias,
Universidade do Porto, Portugal}\\
{\small\texttt{jpp@fc.up.pt}}
}

\graphicspath{{figures_flow/}}

\begin{document}

\trtitle{\ptitleA}
\trauthor{\pauthor}
\trnumber{\pnumber}
\mkcoverpage

\title{\textbf{\ptitleB}}
\author{\pauthor}
\date{26 September, 2013}
\maketitle

\begin{abstract}
The cutting stock problem with binary patterns (0-1 CSP)
is a variant of CSP that usually appears as a relaxation 
of 2D and 3D packing problems.
We present an exact method, based on an arc-flow formulation
with side constraints, for solving 0-1 CSP
by simply representing all the patterns in a very compact graph.

Gilmore-Gomory's column generation approach is usually
used to compute strong lower bounds for 0-1 CSP.
We report a computational comparison between the arc-flow approach 
and the Gilmore-Gomory's approach.
\end{abstract}


\section{Introduction}

The cutting stock problem (CSP) is a combinatorial NP-hard problem 
(see, e.g.,~\citealt{Garey:1979:CIG:578533}) in which pieces of different
widths must be cut from rolls in such a way that the waste is minimized. 
In this problem, we are given the width $W$ of the rolls,
a set of $m$ piece widths $w$ and their demands $b$.
Since all the rolls have the same width
and demands must be met, the objective is equivalent to 
minimizing the number of rolls that are used.
In the cutting stock problem with binary patterns (0-1 CSP) items
of each type may be cut at most once in each roll.
In the 0-1 CSP, pieces are identified by their types and some types may have the
same width.

The 2D strip packing problem (SPP) is another combinatorial NP-hard problem.
In this problem, we are given a half-open strip of width $W$
and we want to minimize the height needed to pack a set of 2D rectangular items.
An known lower bound for SPP is the bar relaxation
of~\cite{Scheithauer1999}, 
which corresponds to the minimum number of one-dimensional packing patterns in width 
direction (bar patterns with length one) that are required to hold all the items. 
This relaxation corresponds to a cutting stock problem with binary patterns (0-1 CSP)
since each bar contains at most one slice of each item.
For 3D problems, there are bar and slice relaxations (see, e.g.,~\citealt{OPPBounds}), 
both of which are also 0-1 CSP problems.

The 0-1 CSP can be reduced and solved as a vector packing problem with $m+1$ 
dimensions, one for the capacity constraint and $m$ binary dimensions to ensure that 
each pattern contains at most one item of each type.
\cite{MThesisBrandao} presents a graph compression method
for vector packing arc-flow graphs that usually leads to large reductions in the graph size.
In this paper, we use a different method for introducing binary 
constraints in arc-flow models that only requires a single additional dimension.
The main contributions of this paper are: 
we present a exact method for solving 0-1 CSP with binary patterns
which can be easily generalized for vector packing with binary patterns.
This generalization allows, for instance, modeling 0-1 CSP with conflicts,
which is another problem that usually appears when solving 2D and 3D packing problems.

The remainder of this paper is organized as follows.
Section~\ref{sec:models} presents underlying mathematical optimization models,
and Section~\ref{sec:method} the arc-flow approach to 0-1 CSP.
Some computational results are presented in Section~\ref{sec:results}
and Section~\ref{sec:conclusions} presents the conclusions.

\section{Mathematical optimization models}
\label{sec:models}

\cite{gomory1} proposed the following model for the standard CSP.
A combination of orders in the width of the
roll is called a cutting pattern.  Let column vectors $a^j = (a_1^j,
\ldots, a_m^j)^{\top}$ represent all possible cutting patterns~$j$. The
element $a_i^j$ represents the number of rolls of width $w_i$ obtained
in cutting pattern $j$.  
Note that in the standard CSP, the patterns are non-negative integer vectors;
in the 0-1 CSP, only binary patterns are allowed.
Let $x_j$ be a decision variable that
designates the number of rolls to be cut according to cutting pattern
$j$.  The 0-1 CSP can be modeled in terms of these
variables as follows:
\begin{alignat}{3}
  & \mbox{minimize }          && \sum_{j \in J} x_j \label{eq:cuttingstock1}\\  
  & \mbox{subject to } \qquad && \sum_{j \in J} a_i^j x_j \geq b_i, \qquad && i=1,\ldots,m, \label{eq:cuttingstock2}\\
  &                           && x_{j} \geq 0, \mbox{ integer},    \qquad && \forall j \in J, \label{eq:cuttingstock3}
\end{alignat}
where $J$ is the set of valid cutting patterns that satisfy:
\begin{equation}
\sum_{i = 1}^{m} a_i^j  w_i \leq W \mbox{ and }  a_i^j \in \{0,1\}. \label{eq:cuttingstock4}
\end{equation}

It may be impractical to enumerate all the columns in the
previous formulation, as their number may be very large, even for moderately 
sized problems. To tackle this problem, \cite{gomory2} proposed column
generation.

Let the linear optimization of Model~(\ref{eq:cuttingstock1})-(\ref{eq:cuttingstock3})
be the restricted master problem.
At each iteration of the column generation process, a subproblem
is solved and a column (pattern) is introduced in the restricted master problem
if its reduced cost is strictly less than zero.
The subproblem, which is a knapsack problem, is the following:
\begin{alignat}{3}
  & \mbox{minimize}          && 1-\sum_{i = 1}^{m} c_i a_i \label{eq:subprob1}\\  
  & \mbox{subject to } \qquad && \sum_{i = 1}^{m} w_i a_i \leq W\label{eq:subprob2}\\
  &                           && a_i \in \{0,1\},\ i=1,\ldots,m,  \label{eq:subprob3}
\end{alignat}
where: $c_i$ is the shadow price of the demand constraint of item $i$
obtained from the solution of the linear relaxation of the restricted master problem,
and $a = (a_1, \ldots, a_m)$ is a cutting pattern 
whose reduced cost is given by the objective function.

The column generation process for this method can be 
summarized as follows. 
We start with a small set of patterns (columns),
which can be composed by $m$ patterns, each containing a single, different item.
Then we solve the linear relaxation of the restricted master problem with
the current set of columns. At each iteration,
a knapsack subproblem is solved and the pattern $a^*$ from its solution
is introduced in the restricted master problem.
Simplex iterations are then performed to update the solution
of the master problem.
This process is repeated until no pattern with negative reduced cost
is found. At the end of this process, we have the optimal solution of the linear
relaxation of the model~(\ref{eq:cuttingstock1})-(\ref{eq:cuttingstock3}).

Model~(\ref{eq:subprob1})-(\ref{eq:subprob3}) corresponds to a 0-1 knapsack problem
that can be solved in pseudo-polynomial time with dynamic programming
using Algorithm~\ref{alg:01kp}, that runs in $\BigO(Wm)$ time.

\begin{figure}[h!tbp]
\caption{Dynamic programming graph for a knapsack instance with capacity 8 and items of sizes 4, 3, 2.\label{ex:dpgraph}}

  \centering
  \includegraphics[scale=1]{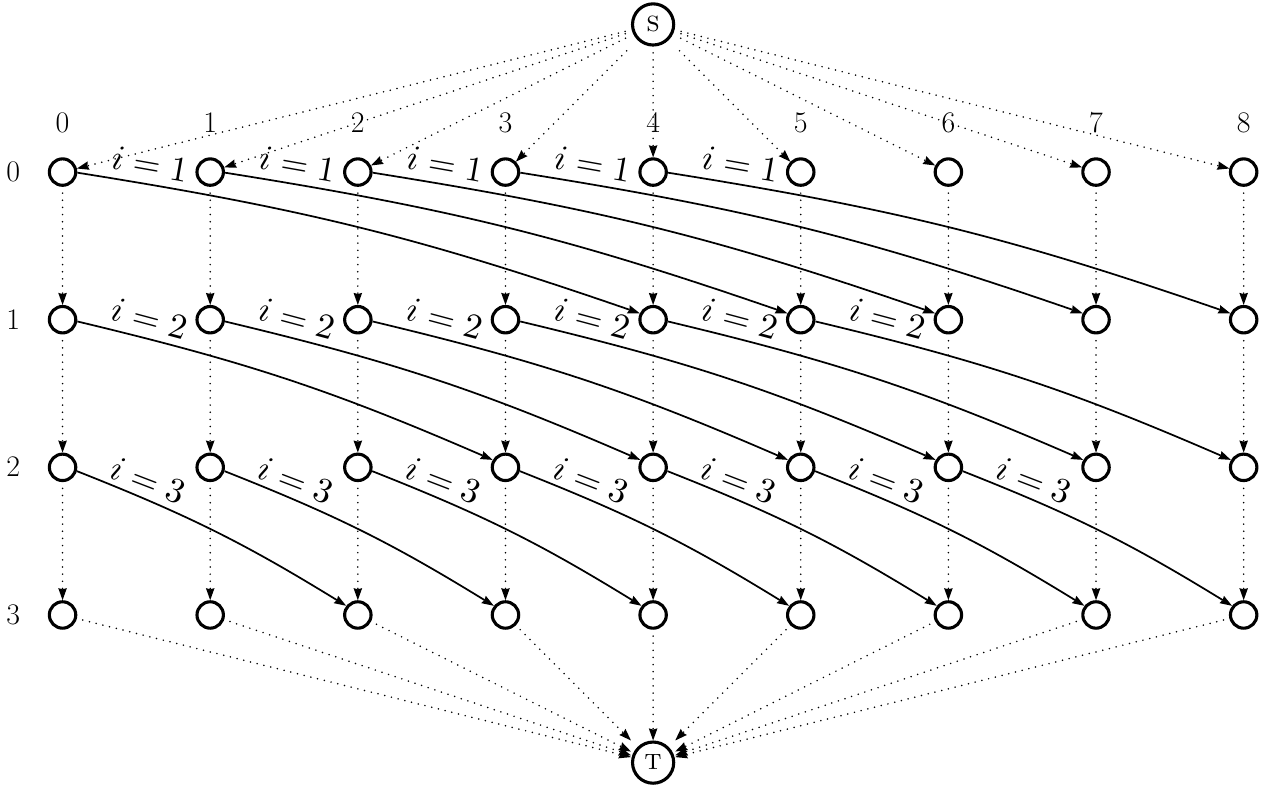}

\caption*{\footnotesize 
}
\end{figure}

Consider a 0-1 CSP instance with bins of capacity $W = 8$ and items of sizes 4, 3, 2. 
The dynamic programming search space of Algorithm~\ref{alg:01kp} is represented 
in Figure~\ref{ex:dpgraph}, and it corresponds to a directed acyclic graph in
which every valid packing pattern is represented as a path from $\vS$ to $\vT$.
\cite{MThesisBrandao} presents a general arc-flow formulation,
equivalent to Model~(\ref{eq:cuttingstock1})-(\ref{eq:cuttingstock4}),
that can be used to solve 0-1 CSP directly as a minimum flow problem 
between $\vS$ and $\vT$, with 
additional constraints enforcing the sum of the flows in the arcs of
each item to be greater than or equal to the corresponding demand.
This general arc-flow formulation is a generalization of the model proposed in \cite{Valerio:01} and 
it only requires a directed acyclic graph 
containing every valid packing pattern represented as a path between two vertices
to solve the corresponding cutting stock problem.
The lower bound provided by this formulation is the same as the one 
provided by Gilmore-Gomory's model with the same set of patterns (see, e.g., \citealt{MThesisBrandao}).
A simplified version of the general arc-flow formulation is the following:
\begin{alignat}{3}
  & \mbox{minimize }   && z  \label{eq:new1}\\  
  & \mbox{subject to } \qquad && \sum_{(u,v,i)\in A} f_{uvi}- \sum_{(v,u^{\prime},i) \in A}f_{vu^{\prime} i} = &&
  \left\{ 
    \begin{array}{rl}
      -z & \mbox{if } v = \vS, \\
      z & \mbox{if } v = \vT,\\
      0 & \mbox{for } v \in V \setminus \{\vS, \vT \},\\
    \end{array} \label{eq:new2}
  \right.\\
      &&         & \sum_{(u, v, i) \in A} f_{uvi} \geq b_i, && i \in \{1,\ldots,m\}, \label{eq:new3}\\ 
      &&         & f_{uvi} \geq 0, \mbox{ integer}, && \forall (u,v,i) \in A, \label{eq:new4}      
\end{alignat}
where: $m$ is the number of different items;
$b_i$ is the demand of the $i$-th item;
$V$ is the set of vertices,
$\vS$ is the source vertex and $\vT$ is the target;
$A$ is the set of arcs,
each arc having components $(u, v, i)$
corresponding to an arc between nodes $u$ and $v$
that contributes to the demand of the $i$-th item;
arcs $(u,v,i=0)$ are loss arcs that represent unoccupied portions of the patterns;
$f_{uvi}$ is the amount of flow along the arc $(u, v, i)$;
and $z$ is a variable that can be seen as a feedback arc from vertex $\vT$ to $\vS$.
Note that this formulation allows multiple arcs between the same pair
of vertices.

\begin{algorithm}[!h]
\caption{0-1 Knapsack Algorithm}\label{alg:01kp}
\SetKwInOut{Input}{input}
\SetKwInOut{Output}{output}

\SetKwFunction{Build}{build}
\SetKwFunction{Knapsack}{knapsack}
\SetKwFunction{Sort}{sort}
\SetKwFunction{Key}{key}
\SetKwFunction{Reversed}{reversed}
\SetKwData{vardp}{dp}
\SetKwData{varx}{x}
\SetKwData{vari}{i}
\SetKwData{varmx}{mx}
\SetKwData{varmi}{mi}
\SetKwData{NIL}{NIL}

\SetKwBlock{Function}{}{}

\Input{$m$ - number of different items; $w$ - item sizes; $v$ - item values; $W$ - capacity limit}
\Output{maximum profit}

\textbf{function} $\Knapsack(m, w, v, W)$:
\Function{
\lFor(\tcp*[f]{base cases}){$p = 0$ \KwTo $W$}{$\vardp[0, p]\gets 0$} 
\For(\tcp*[f]{for each item}){$i = 1$ \KwTo $m$}{
    \For(\tcp*[f]{for each position}){$p = 0$ \KwTo $W$}{    
        \lIf(){$p < w_i$}{$\vardp[i, p]\gets \vardp[i-1,p] $}
        \lElse(){$\vardp[i, p]\gets \max(\vardp[i-1,p], \vardp[i-1,p-w_i]+v_i)$}
    }
} 
\Return $\max_{p=0..W} \{\vardp[m, p]\}$\;
}
\end{algorithm}

\FloatBarrier

\section{Arc-flow formulation with graph compression for 0-1 CSP}
\label{sec:method}

\begin{figure}[h!tbp]
\caption{Initial graph/Step-1 graph.\label{ex:compression0}}

  \centering
  \includegraphics[scale=1]{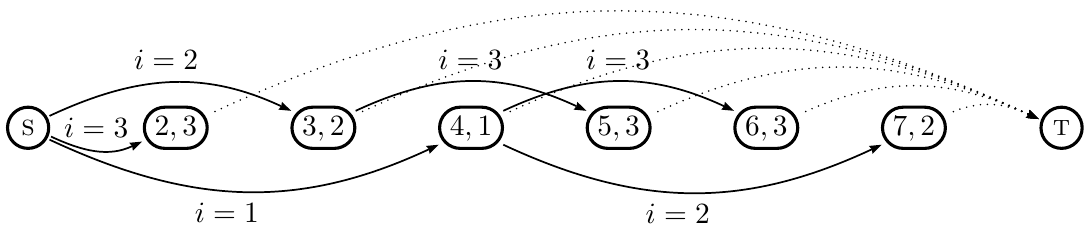}

\caption*{\footnotesize
Graph corresponding to a 0-1 CSP instance
with bins of capacity $W = 8$ and items of sizes 
4, 3, 2 with demands 3, 2, 5, respectively.
}
\end{figure}

Consider a 0-1 CSP instance
with bins of capacity $W = 8$ and items of sizes 
4, 3, 2 with demands 3, 2, 5, respectively.
Since a pattern is a set of items, in order to avoid redundant
patterns we will assume that this set is ordered, and consider
only patterns in which the items are ordered by decreasing values
of the index in this set.
Figure~\ref{ex:compression0} shows a graph which contains every valid packing
pattern (respecting that order) for this instance represented as path from the 
source $\vS$ to the target $\vT$.
In this graph, a node label $(a^{\prime},b^{\prime})$
means that every sub-pattern from the source to the node
uses no more than $a^{\prime}$ space and contains no item with an index higher than $b^{\prime}$.
This graph can be seen as the Step-1 graph of \cite{MThesisBrandao}'s method.
The dashed arcs are the loss arcs that represent unoccupied portions of the patterns,
which can be seen as items with no width ($w_0 = 0$).
In a general arc-flow graph, 
a tuple $(u, v, i)$ corresponds to an arc between nodes 
$u$ and $v$ associated with the $i$-th item.
Note that, for each pair of nodes $(u, v)$, multiple arcs, each associated with a different item,
are allowed.

\cite{MThesisBrandao} presents a three-step graph 
compression method whose first step consists of breaking the symmetry
by dividing the graph into levels.
The way used to represent patterns in Figure~\ref{ex:compression0}
does not allow symmetry. 
The graph division by levels does not improve this, but it
increases the flexibility of the graph
and usually allows substantial improvements in the compression ratio.
Note that this graph division by levels does not exclude any valid packing
pattern that respects the order (see, e.g., \citealt{MThesisBrandao}). 
It is easy to check that, excepting loss,
for every valid packing pattern in the initial graph, 
there is a corresponding path in the Step-2 graph.

\begin{figure}[h!tbp]
\caption{Graph with levels/Step-2 graph.\label{ex:compression1}}

  \centering
  \includegraphics[scale=1]{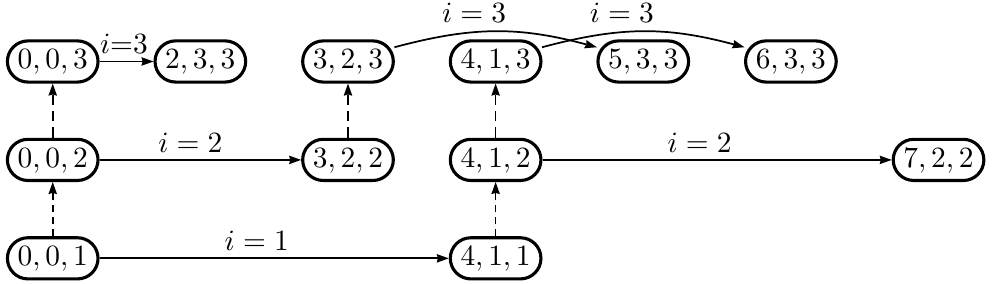}

\caption*{
\footnotesize 
The Step-2 graph has 13 nodes and 22 arcs
(considering also the final loss arcs connecting internal nodes to $\vT$,
which were omitted). In this graph, $\vS=(0,0,1)$.
}
\end{figure}
\FloatBarrier 

\begin{figure}[h!tbp]
\caption{Step-3 graph (after the main compression step).\label{ex:compression2}}

  \centering
  \includegraphics[scale=1]{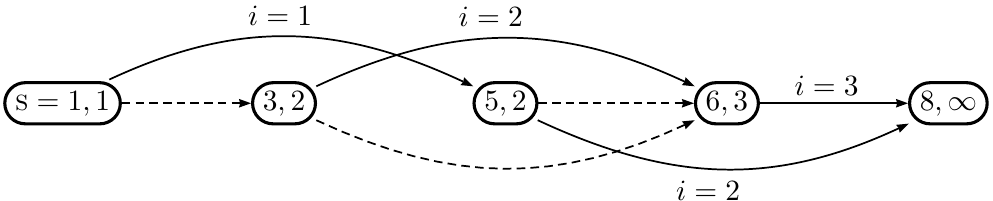}

\caption*{\footnotesize The Step-3 graph has 8 nodes and 17 arcs 
(considering also the final loss arcs connecting internal nodes to~$\vT$,
which were omitted).
}
\end{figure}

In the main compression step, a new graph is constructed using the longest
path to the target in each dimension.
In our case, we use in the first dimension the longest path from the node to the target, 
and in the second the highest item index that appears in any path from
the node to the target.
Let ($\varphi(u)$, $\psi(u)$) be the label of node $u$ in the
first and second dimensions, respectively.
We define $\varphi(u)$ and $\psi(u)$ as follows:
\begin{alignat}{3}
& \varphi(u) && = && \left\{ \begin{array}{ll}
                W & \mbox{if }  u = \vT, \\
                \min_{(u',v,i) \in A: u'=u}\{\varphi(v) - w_i\}  & \mbox{otherwise.}\\
                \end{array}\right.\\
& \psi(u) && = && \left\{ \begin{array}{ll}
                \infty & \mbox{if }  u = \vT, \\
                \min(\min_{(u',v,i) \in A: u'=u}\{\psi(v)\}, \min_{(u,v,i\neq0) \in A}\{i\}) & \mbox{otherwise.}\\
                \end{array}\right.                
\end{alignat}
In the paths from $\vS$ to $\vT$ in Step-2 graph usually there is some float.
In this process, we are moving this float as much as possible to the beginning of the paths.
The label in each dimension of every node $u$ corresponds to the highest value in each dimension
where the sub-patterns from $u$ to $\vT$ can start
so that the constraints are not violated.
Figure~\ref{ex:compression2} shows the graph that results from applying
the main compression step to the graph of Figure~\ref{ex:compression1}.

\FloatBarrier 

\begin{figure}[h!tbp]
\caption{Step-4 graph (after the final compression step).\label{ex:compression3}}

  \centering
  \includegraphics[scale=1]{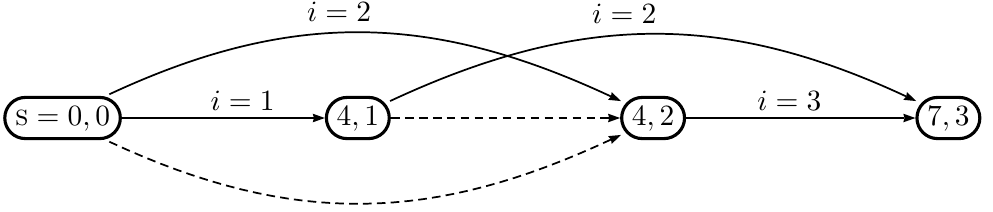}

\caption*{\footnotesize The Step-4 graph has 5 nodes and 9 arcs
(considering also the final loss arcs connecting internal nodes to $\vT$,
which were omitted).
The initial Step-1 graph had 8 nodes and 12 arcs.}
\end{figure}

In the final compression step, a new graph is constructed once more,
using in the first dimension the longest path from the source to the current node, 
and in the second the highest item index that appears in any path from
the source to the current node.
Figure~\ref{ex:compression3} shows the final graph.
Let ($\varphi^{\prime}(v)$, $\psi^{\prime}(v)$) be the label of node $v$ in the
first and second dimensions, respectively.
We define $\varphi^{\prime}(v)$ and $\psi^{\prime}(v)$ as follows:
\begin{alignat}{3}
& \varphi^{\prime}(v) && = && \left\{ \begin{array}{ll}
                0 & \mbox{if }  v = \vS, \\
                \max_{(u,v',i) \in A : v'=v}\{\varphi^{\prime}(u) + w_i\}  & \mbox{otherwise.}\\
                \end{array}\right.\\
& \psi^{\prime}(v) && = && \left\{ \begin{array}{ll}
                0 & \mbox{if }  v = \vS, \\
                \max_{(u,v',i) \in A : v'=v}\{\max(\psi^{\prime}(u),i)\}  & \mbox{otherwise.}\\
                \end{array}\right.                
\end{alignat}

\FloatBarrier 

The problem of minimizing the flow $z$ in the final graph is solved by a general-purpose
mixed-integer optimization solver.
The graph compression usually leads to very small graphs.
However, Step-2 graph can be large and 
the three-step graph compression method may take too long.
Note that 0-1 CSP problems may need to be solved many times, for instance,
as bar or slices relaxations during the solution
of 2D and 3D packing problems, and hence it is very important to build
the models as quick as possible.
Algorithm~\ref{alg:build} builds the Step-3 graph directly
using a recursive procedure with memoization.
The base idea for this algorithm comes from the fact that
in the main compression step the label of any node only depends
on the labels of the two nodes to which it is connected (a node in its level and another in level above).
After directly building the Step-3 graph from the instance's data using this algorithm, we just need 
to apply the final compression step to obtain the final graph. 
This method allows us to obtain arc-flow models even for large benchmark instances in a few milliseconds.

\begin{algorithm}[!h]
\caption{Direct Step-3 Graph Construction Algorithm}\label{alg:build}
\SetKwInOut{Input}{input}
\SetKwInOut{Output}{output}

\SetKwFunction{Build}{build}
\SetKwFunction{BuildGraph}{buildGraph}
\SetKwFunction{Sort}{sort}
\SetKwFunction{Key}{key}
\SetKwFunction{Reversed}{reversed}
\SetKwData{vardp}{dp}
\SetKwData{varx}{x}
\SetKwData{vari}{i}
\SetKwData{varmx}{mx}
\SetKwData{varmi}{mi}
\SetKwData{NIL}{NIL}

\SetKwBlock{Function}{}{}

\Input{$m$ - number of different items; $w$ - item sizes; $W$ - capacity limit}
\Output{$V$ - set of vertices; $A$ - set of arcs; $\vS$ - source; $\vS$ - target}

\textbf{function} $\BuildGraph(m, w, W)$:
\Function{
$V \gets \{\ \}$; $A \gets \{\ \}$\;
$\vardp[x', i'] \gets \NIL, \mbox{ \bf for all } x', i'$\;
\textbf{function} $\Build(x, i)$:
\Function{
\If(\tcp*[f]{avoid repeating work }){$\vardp[x,i] \neq \NIL$}{
    \Return $\vardp[x,i]$\;
}
$(u_x, u_i) \gets (W, m+1)$\;
\If(\tcp*[f]{option 1: do not use the current item}){$i < m$}{
    $(up_x, up_i) \gets \Build(x, i+1)$\;
    $(u_x, u_i) \gets (up_x, up_i)$\;
}
\If(\tcp*[f]{option 2: use the current item}){$x+w_i \leq W$}{
    $(v_x, v_i) \gets \Build(x+w_i, i+1)$\;
    $(u_x, u_i) \gets (\min(u_x, v_x-w_i), i)$\;
    $A \gets A \cup \{((u_x, u_i), (v_x, v_i),i)\}$\;
    $V \gets V \cup \{(u_x, u_i), (v_x, v_i)\}$\;
    \If(){$i < m${ \bf and }$u \neq up$}{
        $A \gets A \cup \{((u_x, u_i),(up_x, up_i),0)\}$\;
        $V \gets V \cup \{(up_x, up_i)\}$\;        
    }
}
$\vardp[x,i] \gets (u_x, u_i)$\;

\Return $(u_x, u_i)$\;
}
$\vS \gets \Build(0,1)$\tcp*{build the graph}
$V \gets V \cup \{\vT\}$\;
$A \gets A \cup \{(u, \vT, 0)\ |\ u \in V\setminus \{\vS, \vT\} \}$\tcp*{connect internal nodes to the target}
\Return $(G = (V, A), \vS, \vT)$\;
}
\end{algorithm}

\pagebreak
\section{Computational results}
\label{sec:results}

We used the arc-flow formulation to compute the bar relation
of the two-dimensional bin packing test data of~\cite{Lodi1999}.
This data set is composed by ten classes;
the first six were proposed by~\cite{bw-tdfbp-87},
and the last four were proposed by~\cite{Martello1998}.
In the bar relaxations, the item widths are the original widths
and the demand is given by the original height of the item.
The characteristics of the first six classes are summarized as follows:
Class I: $W=10$, $w_i \in [1,10]$ and $b_i \in [1,10]$;
Class II: $W=30$, $w_i \in [1,10]$ and $b_i \in [1,10]$;
Class III: $W=40$, $w_i \in [1,35]$ and $b_i \in [1,35]$;
Class IV: $W=100$, $w_i \in [1,35]$ and $b_i \in [1,35]$;
Class V: $W=100$, $w_i \in [1,100]$ and $b_i \in [1,100]$;
Class VI: $W=300$, $w_i \in [1,100]$ and $b_i \in [1,100]$.
In each of the previous classes, all the item sizes are generated in the same interval.
In the next four classes, item sizes are generated from different intervals.
The characteristics of 70\% of the items in each class are summarized as follows:
Class VII: $W=100$, $w_i \in [66,100]$ and $b_i \in [1,50]$;
Class VIII: $W=100$, $w_i \in [1,50]$ $b_i \in [66,100]$;
Class IX: $W=100$, $w_i \in [50,100]$ $b_i \in [50,100]$;
Class X: $W=100$, $w_i \in [1,50]$ $b_i \in [1,50]$.

Table~\ref{tab:results} presents a comparison
between Gilmore-Gomory's column generation approach and the
arc-flow formulation approach.
The meaning of each column is as follows: $W$ - strip width;
$m$ - number of different items;
CSP - optimal standard CSP solution;
$\zlp$ - LP lower bound;
$\zip$ - optimal 0-1~CSP solution;
$\ncols$ -  number of columns used in the column generation approach;
$\#v$, $\#a$ - number of vertices and arcs in the arc-flow graph;
$\%v$, $\%a$ - ratios between the numbers of vertices and arcs of the arc-flow
graph and the ones of the straightforward dynamic programming graph presented in Section~\ref{sec:models};
$\tgg$ - column generation run time;
$\tlp$, $\tip$ - run time of the arc-flow approach (linear relaxation and integer solution);
$\bestgg$, $\bestlp$ - number of instances that were solved in $\min(\tgg, \tlp)$ seconds
using each method.
The values shown are averages over the 10 instances in each class. 
CPU times were obtained using a computer with two Quad-Core Intel 
Xeon at 2.66GHz, running Mac OS X 10.8.0.
All the algorithms were implemented in \texttt{C++},
and \texttt{Gurobi} 5.0.0 (\citealt{gurobi}), a state-of-the-art mixed integer programming solver,
was used as LP and ILP solver.
The source code is available online\footnote{\url{http://www.dcc.fc.up.pt/\~fdabrandao/code}}.

\setlength{\tabcolsep}{2pt}
\begin{table}[h!tbp]
\centering\begin{threeparttable}[b]
\scriptsize 
\caption{Arc-flow results.\label{tab:results}}
\begin{tabular}{lrrrrrrrrrrrrrrr}
\hline
\up\down
Class & $W$ & $m$ & CSP & $\zip$ & $\zlp$ & $\ncols$ & $\#v$ & $\#a$ & $\%v$ & $\%a$ & $\tgg$ & $\tlp$ & $\tip$ & $\bestgg$ & $\bestlp$\\
\hline	
\up
I & 10 & 20 & 60.4 & 60.8 & 60.6 & 47.8 & 23.6 & 87.3 & 0.1 & 0.24 & 0.01 & 0.01 & 0.01 & 0 & 10\\
I & 10 & 40 & 121.2 & 121.7 & 121.4 & 105.3 & 56.1 & 214 & 0.12 & 0.31 & 0.02 & 0.01 & 0.01 & 0 & 10\\
I & 10 & 60 & 188.5 & 188.5 & 188.3 & 147.3 & 80.8 & 318.8 & 0.12 & 0.31 & 0.03 & 0.01 & 0.02 & 0 & 10\\
I & 10 & 80 & 262.6 & 262.6 & 262.2 & 189.8 & 105.3 & 417.9 & 0.12 & 0.31 & 0.04 & 0.01 & 0.02 & 0 & 10\\
I & 10 & 100 & 304.8 & 304.8 & 304.3 & 297 & 163.3 & 629.2 & 0.15 & 0.38 & 0.1 & 0.01 & 0.03 & 0 & 10\\
\up
II & 30 & 20 & 19.7 & 19.7 & 19.3 & 45.5 & 131.9 & 436.3 & 0.2 & 0.36 & 0.01 & 0.01 & 0.06 & 9 & 4\\
II & 30 & 40 & 39.1 & 39.1 & 38.7 & 93.1 & 357.9 & 1190.4 & 0.28 & 0.51 & 0.02 & 0.03 & 0.46 & 9 & 1\\
II & 30 & 60 & 60.1 & 60.1 & 59.8 & 139 & 602.8 & 2008.9 & 0.32 & 0.58 & 0.04 & 0.05 & 2.01 & 8 & 2\\
II & 30 & 80 & 83.2 & 83.2 & 82.7 & 188 & 812.9 & 2719.7 & 0.32 & 0.6 & 0.07 & 0.07 & 3.2 & 3 & 7\\
II & 30 & 100 & 100.5 & 100.5 & 99.9 & 240.6 & 1131.1 & 3727.4 & 0.36 & 0.65 & 0.13 & 0.11 & 3.03 & 0 & 10\\
\up
III & 40 & 20 & 160.1 & 161 & 160.8 & 50.4 & 44.3 & 166.7 & 0.05 & 0.12 & 0.01 & 0.01 & 0.01 & 1 & 10\\
III & 40 & 40 & 330.4 & 331.6 & 331.3 & 108.4 & 141.5 & 539.8 & 0.08 & 0.2 & 0.02 & 0.01 & 0.02 & 0 & 10\\
III & 40 & 60 & 502.6 & 503.4 & 503 & 169.3 & 224.5 & 888.9 & 0.09 & 0.23 & 0.04 & 0.02 & 0.05 & 0 & 10\\
III & 40 & 80 & 707.2 & 707.2 & 706.9 & 228.1 & 308.8 & 1215.4 & 0.09 & 0.24 & 0.06 & 0.03 & 0.06 & 0 & 10\\
III & 40 & 100 & 834.8 & 834.9 & 834.4 & 307.1 & 525.4 & 1973.5 & 0.13 & 0.3 & 0.12 & 0.05 & 0.19 & 0 & 10\\
\up
IV & 100 & 20 & 61.4 & 61.4 & 60.9 & 57 & 249.4 & 859.8 & 0.12 & 0.22 & 0.01 & 0.02 & 0.09 & 10 & 0\\
IV & 100 & 40 & 123.9 & 123.9 & 123.3 & 107.1 & 866.7 & 2962.7 & 0.21 & 0.39 & 0.03 & 0.07 & 0.77 & 10 & 0\\
IV & 100 & 60 & 193 & 193 & 192.4 & 158.5 & 1618.3 & 5502.5 & 0.26 & 0.49 & 0.06 & 0.16 & 2.11 & 10 & 0\\
IV & 100 & 80 & 267.2 & 267.2 & 266.7 & 219.5 & 2209.2 & 7528.3 & 0.27 & 0.51 & 0.11 & 0.29 & 3.64 & 10 & 0\\
IV & 100 & 100 & 322 & 322 & 321.4 & 280.2 & 3157 & 10545.7 & 0.31 & 0.57 & 0.2 & 0.46 & 11.51 & 10 & 0\\
\up
V & 100 & 20 & 515.5 & 530.2 & 530 & 42.5 & 39.8 & 150.4 & 0.02 & 0.05 & 0.01 & 0.01 & 0.01 & 2 & 10\\
V & 100 & 40 & 1061.2 & 1069.1 & 1069.1 & 89.5 & 164.8 & 648.8 & 0.04 & 0.1 & 0.02 & 0.01 & 0.02 & 2 & 9\\
V & 100 & 60 & 1625.1 & 1635.4 & 1635.2 & 153 & 307.8 & 1246.2 & 0.05 & 0.13 & 0.04 & 0.03 & 0.04 & 1 & 9\\
V & 100 & 80 & 2281.4 & 2283.8 & 2283.4 & 197.9 & 461.6 & 1839.7 & 0.06 & 0.15 & 0.05 & 0.03 & 0.06 & 0 & 10\\
V & 100 & 100 & 2634.3 & 2634.8 & 2634.4 & 296.1 & 875.8 & 3320.6 & 0.09 & 0.22 & 0.11 & 0.07 & 0.17 & 0 & 10\\
\up
VI & 300 & 20 & 159.9 & 159.9 & 159.4 & 69.2 & 412.6 & 1464.4 & 0.07 & 0.13 & 0.02 & 0.03 & 0.19 & 10 & 0\\
VI & 300 & 40 & 323.5 & 323.5 & 323.1 & 126 & 2046.3 & 7010.3 & 0.17 & 0.31 & 0.04 & 0.24 & 2.2 & 10 & 0\\
VI & 300 & 60 & 505.1 & 505.1 & 504.7 & 186.6 & 4220.4 & 14361.9 & 0.23 & 0.43 & 0.09 & 0.64 & 14.68 & 10 & 0\\
VI & 300 & 80 & 699.7 & 699.7 & 699.1 & 255.4 & 5908.6 & 20153.3 & 0.24 & 0.45 & 0.17 & 1.1 & 20.45 & 10 & 0\\
VI & 300 & 100 & 843.8 & 843.8 & 843.1 & 329.1 & 8646.9 & 28900.2 & 0.28 & 0.52 & 0.33 & 2.34 & 75.14 & 10 & 0\\
\up
VII & 100 & 20 & 484.9 & 500.6 & 500.6 & 25.9 & 9.2 & 37.5 & 0 & 0.01 & 0.01 & 0.01 & 0.01 & 5 & 10\\
VII & 100 & 40 & 1050.7 & 1054.7 & 1054.6 & 55.2 & 21.2 & 98.6 & 0.01 & 0.02 & 0.01 & 0.01 & 0.01 & 3 & 10\\
VII & 100 & 60 & 1513.3 & 1525 & 1524.9 & 87.9 & 48 & 219.9 & 0.01 & 0.03 & 0.01 & 0.01 & 0.01 & 2 & 10\\
VII & 100 & 80 & 2210.6 & 2220.8 & 2220.7 & 117.4 & 78.1 & 374.4 & 0.01 & 0.03 & 0.01 & 0.01 & 0.01 & 0 & 10\\
VII & 100 & 100 & 2633.3 & 2642.7 & 2642.6 & 151.1 & 152.8 & 656.5 & 0.01 & 0.05 & 0.02 & 0.01 & 0.02 & 0 & 10\\
\up
VIII & 100 & 20 & 437.4 & 437.9 & 437.3 & 60.6 & 97.7 & 364.9 & 0.05 & 0.1 & 0.01 & 0.01 & 0.02 & 2 & 9\\
VIII & 100 & 40 & 922.7 & 924.7 & 924.1 & 121.2 & 351.4 & 1325.4 & 0.08 & 0.19 & 0.03 & 0.03 & 0.12 & 6 & 5\\
VIII & 100 & 60 & 1360.9 & 1360.9 & 1360.5 & 191.8 & 776.4 & 2832.3 & 0.13 & 0.28 & 0.06 & 0.08 & 0.56 & 9 & 2\\
VIII & 100 & 80 & 1911.4 & 1911.4 & 1911.1 & 257.5 & 1081.3 & 3977.6 & 0.13 & 0.3 & 0.11 & 0.14 & 0.9 & 8 & 2\\
VIII & 100 & 100 & 2363 & 2363 & 2362.5 & 330.5 & 1708.7 & 6085.3 & 0.17 & 0.36 & 0.18 & 0.26 & 2.65 & 10 & 0\\
\up
IX & 100 & 20 & 1103.6 & 1106.8 & 1106.8 & 27.5 & 10.3 & 42.4 & 0 & 0.01 & 0.01 & 0.01 & 0.01 & 5 & 10\\
IX & 100 & 40 & 2180.9 & 2190 & 2189.9 & 53.3 & 22.9 & 110.6 & 0.01 & 0.02 & 0.01 & 0.01 & 0.01 & 2 & 10\\
IX & 100 & 60 & 3394 & 3410.4 & 3410.4 & 79.8 & 40.4 & 202.1 & 0.01 & 0.02 & 0.01 & 0.01 & 0.01 & 3 & 9\\
IX & 100 & 80 & 4563.9 & 4584.4 & 4584.4 & 106.5 & 82.3 & 402.2 & 0.01 & 0.04 & 0.01 & 0.01 & 0.01 & 2 & 9\\
IX & 100 & 100 & 5415.6 & 5431.2 & 5431.2 & 146.6 & 190.9 & 851.4 & 0.02 & 0.06 & 0.02 & 0.01 & 0.02 & 1 & 10\\
\up
X & 100 & 20 & 341.4 & 347.1 & 346.9 & 51.9 & 86.3 & 323.1 & 0.04 & 0.09 & 0.01 & 0.01 & 0.01 & 1 & 10\\
X & 100 & 40 & 653.7 & 654.6 & 654.2 & 121.4 & 315.5 & 1191.8 & 0.08 & 0.18 & 0.03 & 0.03 & 0.07 & 2 & 8\\
X & 100 & 60 & 918.4 & 919.5 & 919.2 & 213.3 & 752.9 & 2741.7 & 0.12 & 0.27 & 0.06 & 0.08 & 0.23 & 7 & 3\\
X & 100 & 80 & 1185.3 & 1186 & 1185.3 & 311.5 & 1084.7 & 3970.9 & 0.13 & 0.3 & 0.13 & 0.16 & 0.39 & 5 & 5\\
\down
X & 100 & 100 & 1480.7 & 1480.7 & 1480.4 & 394.9 & 1650.1 & 5931.2 & 0.16 & 0.35 & 0.22 & 0.32 & 2.14 & 8 & 2\\
\hline	
\up\down																 											
&\multicolumn{2}{c}{Summary} & 1030.4 & 1033.7 & 1033.3 & 156.6 & 889.1 & 3095.4 & 0.12 & 0.25 & 0.06 & 0.14 & 2.95 & 216 & 316\\
\hline																	 											
\end{tabular}
\end{threeparttable}
\end{table}

The arc-flow formulation is faster than Gilmore-Gomory's approach in 284 instances,
both methods present the same run time in 32 instances
and Gilmore-Gomory's approach is faster than the arc-flow approach in 184 instances.
The hardest instances for the arc-flow approach are the ones of classes IV and VI
since they mostly contain small items that lead to very long patterns.
Graphs associated to this type of instances tend to be larger
and the run times tend to be higher due to the larger
number of variables and constraints.

In Gilmore-Gomory's approach, we used dynamic programming to solve 
the knapsack problems. The arc-flow graph can also be used to solve knapsack
problems. However, since $W$ and $m$ are small in this data set, 
a straightforward dynamic programming solution is faster
due its very low constant factors and its good caching behavior.
For instances with larger values of $W$ and $m$, the use of the 
arc-flow graph to solve the knapsack problems may improve substantially 
the run time.

The graph compression achieved very good reductions in these instances.
The average ratios between the numbers of vertices and arcs of the arc-flow
graph and the ones of the straightforward dynamic programming graph
are $25\%$ and $6\%$, respectively. Moreover, there are instances
whose corresponding compressed graph is more than one hundred times smaller
than the dynamic programming graph.
There are ten instances whose corresponding graph became so small
that it was faster to solve them exactly than 
to compute the linear relaxation using Gilmore-Gomory's approach.
All the instances from this data set were solved exactly 
using the arc-flow method in less than 3 seconds, on average.

The standard CSP can be used as a relaxation of 0-1 CSP.
However, the bound provided by 0-1 CSP is much stronger than the one provided by standard CSP.
In this data set, the maximum difference between the two bounds is 49,
a clear indication of the superiority of 0-1 CSP is better for assessing bar and slice relaxations.

There have been many studies (see, e.g., \citealt{Scheithauer95themodified}, \citealt{Scheithauer1997}) 
about the integrality gap of Gilmore-Gomory's model for the standard cutting stock problem.
The largest gap known so far is $7/6$ and it was found by~\cite{Rietz2002_bounds}.
\cite{Scheithauer1997} conjecture that
$\zip \leq \lceil \zlp \rceil + 1$ for the standard CSP.
For the 0-1 CSP, Gilmore-Gomory model's bounds are also very strong.
Note that our arc-flow model is equivalent to Gilmore-Gomory's model
and hence the lower bounds provided by the linear relaxations are the same.
The largest gap between the linear relaxation of the 0-1 CSP formulations
and the exact solution is 0.99 in this data set.
Therefore, the linear relaxation is usually enough 
when 0-1 CSP is used as a relaxation of other problems.
Nevertheless, for other applications and other problems with binary constraints, 
the exact solution may be very important.

\cite{Kantorovich:1} introduced an assignment-based mathematical programming formulation for CSP,
which can be easily modified for 0-1 CSP by restricting the variables to binary values.
Assignment-based formulations are usually highly symmetric and provide very weak lower bounds 
(see, e.g., \citealt{MThesisBrandao}).
Therefore, this type of models are usually very inefficient in practice.
Using this assignment-based model, we were only able to solve 126 out of the 500 instances within a 10 minute
time limit.
Using the arc-flow formulation, all the instances were easily solved, spending 3 seconds per instance on average,
and none of the instances took longer than 200 seconds to be solved exactly.

\FloatBarrier
\section{Conclusions}
\label{sec:conclusions}

We propose an arc-flow formulation with graph compression
for solving cutting stock problems with binary patterns (0-1 CSP),
which usually appears as a relaxation of 2D and 3D packing problems.
Column generation is usually used to obtain a strong lower bound for
0-1 CSP. 
To the best of our knowledge, we present for the first time
an effective exact method for this problem.
We also report a computational comparison between the arc-flow approach 
and the Gilmore-Gomory's column generation approach for computing lower bounds.

Our method can be easily generalized for vector packing with binary
patterns by using multiple capacity dimensions instead of just one.
This generalization allows, for instance, modeling 0-1 CSP with conflicts,
which is another problem that usually appears when solving 2D and 3D packing problems.

\bibliographystyle{apalike} 
\bibliography{paper}

\begin{thebibliography}{}

\bibitem[Belov et~al., 2009]{OPPBounds}
Belov, G., Kartak, V., Rohling, H., and Scheithauer, G. (2009).
\newblock {One-dimensional relaxations and LP bounds for orthogonal packing}.
\newblock {\em International Transactions in Operational Research},
  16(6):745--766.

\bibitem[Berkey and Wang, 1987]{bw-tdfbp-87}
Berkey, J. and Wang, P. (1987).
\newblock {Two-dimensional finite bin-packing algorithms}.
\newblock {\em J. Oper. Res. Soc.}, 38:423--429.

\bibitem[Brand\~ao, 2012]{MThesisBrandao}
Brand\~ao, F. (2012).
\newblock {Bin Packing and Related Problems: Pattern-Based Approaches}.
\newblock Master's thesis, Faculdade de Ci\^encias da Universidade do Porto,
  Portugal.

\bibitem[Garey and Johnson, 1979]{Garey:1979:CIG:578533}
Garey, M.~R. and Johnson, D.~S. (1979).
\newblock {\em {Computers and Intractability: A Guide to the Theory of
  NP-Completeness}}.
\newblock W. H. Freeman \& Co., New York, NY, USA.

\bibitem[Gilmore and Gomory, 1963]{gomory2}
Gilmore, P. and Gomory, R. (1963).
\newblock {A linear programming approach to the cutting stock problem--part
  {II}}.
\newblock {\em Operations Research}, 11:863--888.

\bibitem[Gilmore and Gomory, 1961]{gomory1}
Gilmore, P.~C. and Gomory, R.~E. (1961).
\newblock {A Linear Programming Approach to the Cutting-Stock Problem}.
\newblock {\em Operations Research}, 9:849--859.

\bibitem[Gu et~al., 2012]{gurobi}
Gu, Z., Rothberg, E., and Bixby, R. (2012).
\newblock {Gurobi Optimizer, Version 5.0.0}.
\newblock (Software program).

\bibitem[Kantorovich, 1960]{Kantorovich:1}
Kantorovich, L.~V. (1960).
\newblock {Mathematical methods of organising and planning production}.
\newblock {\em Management Science}, 6(4):366--422.

\bibitem[Lodi et~al., 1999]{Lodi1999}
Lodi, A., Martello, S., and Vigo, D. (1999).
\newblock {Heuristic and metaheuristic approaches for a class of
  two-dimensional bin packing problems}.
\newblock {\em {INFORMS} Journal on Computing}, 11(4):345--357.

\bibitem[Martello and Vigo, 1998]{Martello1998}
Martello, S. and Vigo, D. (1998).
\newblock {Exact Solution of the Two-Dimensional Finite Bin Packing Problem}.
\newblock {\em Management Science}, 44(3):388--399.

\bibitem[Rietz et~al., 2002]{Rietz2002_bounds}
Rietz, J., Scheithauer, G., and Terno, J. (2002).
\newblock {Tighter bounds for the gap and non-{IRUP} constructions in the
  one-dimensional cutting stock problem}.
\newblock {\em Optimization}, 6:927--963.

\bibitem[Scheithauer, 1999]{Scheithauer1999}
Scheithauer, G. (1999).
\newblock {LP-based bounds for the container and multi-container loading
  problem}.
\newblock {\em International Transactions in Operational Research},
  6(2):199--213.

\bibitem[Scheithauer and Terno, 1995]{Scheithauer95themodified}
Scheithauer, G. and Terno, J. (1995).
\newblock {The Modified Integer Round-Up Property of the One-Dimensional
  Cutting Stock Problem}.

\bibitem[Scheithauer and Terno, 1997]{Scheithauer1997}
Scheithauer, G. and Terno, J. (1997).
\newblock {Theoretical investigations on the modified integer round-up property
  for the one-dimensional cutting stock problem}.
\newblock {\em Operations Research Letters}, 20:93--100.

\bibitem[Val\'erio~de Carvalho, 1999]{Valerio:01}
Val\'erio~de Carvalho, J.~M. (1999).
\newblock {Exact solution of bin-packing problems using column generation and
  branch-and-bound}.
\newblock {\em Ann. Oper. Res.}, 86:629--659.

\end{thebibliography}
 
\end{document}